\pdfoutput=1
\RequirePackage{ifpdf}
\ifpdf 
\documentclass[pdftex]{sigma}
\else
\documentclass{sigma}
\fi

\usepackage{enumerate}
\usepackage[all,2cell]{xy} \UseAllTwocells \SilentMatrices
\numberwithin{equation}{section}
\newtheorem{Theorem}{Theorem}[section]
\newtheorem{Lemma}[Theorem]{Lemma}
{\theoremstyle{definition}
\newtheorem{Definition}[Theorem]{Definition}
\newtheorem{Example}[Theorem]{Example}
}
\DeclareMathOperator{\End}{End}
\DeclareMathOperator{\Sp}{Sp}
\DeclareMathOperator{\Mp}{Mp}

\begin{document}

\newcommand{\arXivNumber}{1410.5529}

\allowdisplaybreaks

\renewcommand{\PaperNumber}{025}

\FirstPageHeading

\ShortArticleName{Metaplectic-c Quantomorphisms}

\ArticleName{Metaplectic-c Quantomorphisms}

\Author{Jennifer VAUGHAN}

\AuthorNameForHeading{J.~Vaughan}

\Address{Department of Mathematics, University of Toronto, Canada}
\Email{\href{mailto:jennifer.vaughan@mail.utoronto.ca}{jennifer.vaughan@mail.utoronto.ca}}

\ArticleDates{Received December 13, 2014, in f\/inal form March 16, 2015; Published online March 24, 2015}

\Abstract{In the classical Kostant--Souriau prequantization procedure, the Poisson algebra of a~symplectic manifold
$(M,\omega)$ is realized as the space of inf\/initesimal quantomorphisms of the prequantization circle bundle.
Robinson and Rawnsley developed an alternative to the Kostant--Souriau quantization process in which the prequantization
circle bundle and metaplectic structure for $(M,\omega)$ are replaced by a~metaplectic-c prequantization.
They proved that metaplectic-c quantization can be applied to a~larger class of manifolds than the classical recipe.
This paper presents a~def\/inition for a~metaplectic-c quantomorphism, which is a~dif\/feomorphism of metaplectic-c
prequantizations that preserves all of their structures.
Since the structure of a~metaplectic-c prequantization is more complicated than that of a~circle bundle, we f\/ind that
the def\/inition must include an extra condition that does not have an analogue in the Kostant--Souriau case.
We then def\/ine an inf\/initesimal quantomorphism to be a~vector f\/ield whose f\/low consists of metaplectic-c
quantomorphisms, and prove that the space of inf\/initesimal metaplectic-c quantomorphisms exhibits all of the same
properties that are seen for the inf\/initesimal quantomorphisms of a~prequantization circle bundle.
In particular, this space is isomorphic to the Poisson algebra $C^\infty(M)$.}


\Classification{53D50; 81S10}

\section{Introduction}

Recall that a~prequantization circle bundle for a~symplectic manifold $(M,\omega)$ consists of a~circle bundle
$Y\rightarrow M$ and a~connection one-form~$\gamma$ on~$Y$ such that $d\gamma=\frac{1}{i\hbar}\omega$.
The Kostant--Souriau quantization recipe with half-form correction requires a~prequantization circle bundle and a~choice
of metaplectic structure for $(M,\omega)$.

Souriau~\cite{sou} def\/ined a~\emph{quantomorphism} between two prequantization circle bundles
$(Y_1,\gamma_1)\rightarrow(M_1,\omega_1)$ and $(Y_2,\gamma_2)\rightarrow(M_2,\omega_2)$ to be a~dif\/feomorphism
$K:Y_1\rightarrow Y_2$ such that $K^*\gamma_2=\gamma_1$.
This condition implies that~$K$ is equivariant with respect to the principal circle actions.
Souriau then def\/ined the \emph{infinitesimal quantomorphisms} of a~prequantization circle bundle $(Y,\gamma)$ to be the
vector f\/ields on~$Y$ whose f\/lows are quantomorphisms.
Kostant~\cite{k1} proved that the space of inf\/initesimal quantomorphisms, which we denote $\mathcal{Q}(Y,\gamma)$, is
isomorphic to the Poisson algebra~$C^\infty(M)$.

The metaplectic-c group is a~circle extension of the symplectic group.
Metaplectic-c quantization, which was developed by Robinson and Rawnsley~\cite{rr1}, is a~variant of Kostant--Souriau
quantization in which the prequantization bundle and metaplectic structure are replaced by a~metaplectic-c structure
$(P,\Sigma)$ and a~prequantization one-form~$\gamma$.
Robinson and Rawnsley proved that metaplectic-c quantization can be applied to all systems that admit metaplectic
quantizations, and to some where the Kostant--Souriau process fails.

In Section~\ref{sec:ks}, we present an explicit construction of the isomorphism from $\mathcal{Q}(Y,\gamma)$ to~$C^\infty(M)$.
In Section~\ref{sec:mpc}, after describing the metaplectic-c prequantization $(P,\Sigma,\gamma)$, we def\/ine
a~me\-ta\-plec\-tic-c quantomorphism, which is a~dif\/feomorphism of metaplectic-c prequantizations that preserves all of their
structures.
Our def\/inition is based on Souriau's, but includes a~condition that is unique to the metaplectic-c context.
We then use the metaplectic-c quantomorphisms to def\/i\-ne~$\mathcal{Q}(P,\Sigma,\gamma)$, the space of inf\/initesimal
metaplectic-c quantomorphisms of $(P,\Sigma,\gamma)$.
We show that every property that was proved for $\mathcal{Q}(Y,\gamma)$ has a~parallel for
$\mathcal{Q}(P,\Sigma,\gamma)$.
In particular, $\mathcal{Q}(P,\Sigma,\gamma)$ is isomorphic to the Poisson algebra $C^\infty(M)$.
The construction in Section~\ref{sec:ks} is used as a~model for the proofs in Section~\ref{sec:mpc}.
We indicate when the calculations are analogous, and when the metaplectic-c case requires additional steps.

Some global remarks concerning notation: for any vector f\/ield $\xi$, the Lie derivative with respect to $\xi$ is written
$L_{\xi}$.
The space of smooth vector f\/ields on a~manifold~$P$ is denoted by~$\mathcal{X}(P)$.
Given a~smooth map $F:P\rightarrow M$ and a~vector f\/ield $\xi\in\mathcal{X}(P)$, we write $F_*\xi$ for the pushforward
of $\xi$ if and only if the result is a~well-def\/ined vector f\/ield on~$M$.
If~$P$ is a~bundle over~$M$, $\Gamma(P)$ denotes the space of smooth sections of~$P$, where the base is always taken to
be the symplectic manifold~$M$.
Planck's constant will only appear in the form~$\hbar$.

\section{Kostant--Souriau quantomorphisms}\label{sec:ks}

In this section, after reviewing the Kostant--Souriau prequantization of a~symplectic mani\-fold~$(M,\omega)$, we construct
a~Lie algebra isomorphism from $C^\infty(M)$ to the space of inf\/initesimal quantomorphisms.
As we have already noted, the fact that these algebras are isomorphic was originally stated by Kostant~\cite{k1} in the
context of line bundles with connection.
His proof can be reconstructed from several propositions across Sections 2--4 of~\cite{k1}.
Kostant's isomorphism is also stated by \'Sniatycki~\cite{sn1}, but much of the proof is left as an exercise.
We are not aware of a~source in the literature for a~self-contained proof that uses the language of principal bundles,
and this is one of our reasons for performing an explicit construction here.

The other goal of this section is to motivate the analogous constructions for a~metaplectic-c prequantization, which
will be the subject of Section~\ref{sec:mpc}.
Each result that we present for Kostant--Souriau prequantization will have a~parallel in the metaplectic-c case.
When the proofs are identical, we will simply refer back to the work shown here, thereby allowing Section~\ref{sec:mpc}
to focus on those features that are unique to metaplectic-c structures.

\subsection{Basic def\/initions and notation}

\subsubsection{Hamiltonian vector f\/ields and the Poisson algebra}

Let $(M,\omega)$ be a~symplectic manifold.
Given $f\in C^\infty(M)$, def\/ine its Hamiltonian vector f\/ield $\xi_f\in\mathcal{X}(M)$~by
\begin{gather*}
\xi_f\lrcorner \omega=df.
\end{gather*}
Def\/ine the Poisson bracket on $C^\infty(M)$~by
\begin{gather*}
\{f,g\}=-\omega(\xi_f,\xi_g),\qquad \forall\, f,g\in C^\infty(M).
\end{gather*}
These choices imply that
\begin{gather*}
\xi_f g=\{f,g\},\qquad \forall\, f,g\in C^\infty(M).
\end{gather*}
A~standard calculation establishes the following fact.

\begin{Lemma}
\label{lem:lemma0}
For all $f,g\in C^\infty(M)$, $[\xi_f,\xi_g]=\xi_{\{f,g\}}$.
\end{Lemma}

\subsubsection{Circle bundles and connections}

Let $Y\stackrel{p}{\longrightarrow}M$ be a~right principal $U(1)$ bundle over a~manifold~$M$.
\begin{itemize}\itemsep=0pt
\item For any $\lambda\in U(1)$, let $R_{\lambda}:Y\rightarrow Y$ represent the right action by~$\lambda$.
That is, $R_{\lambda}(y)=y\cdot\lambda$ for all $y\in Y$.

\item For any $\theta\in\mathfrak{u}(1)$, the Lie algebra of $U(1)$, let $\partial_\theta$ be the vector f\/ield on~$Y$
with f\/low $R_{\exp(t\theta)}$, where $t\in\mathbb{R}$.
In particular, we will consider $\partial_{2\pi i}$.

\end{itemize}
Let~$\gamma$ be a~connection one-form on~$Y$.
By def\/inition,~$\gamma$ is invariant under the right principal action, and for all $\theta\in\mathfrak{u}(1)$,
$\gamma(\partial_\theta)=\theta$.
There is a~two-form~$\varpi$ on~$M$, called the \emph{curvature} of~$\gamma$, such that $d\gamma=p^*\varpi$.

For any $\xi\in\mathcal{X}(M)$, let $\tilde{\xi}$ be the lift of $\xi$ to~$Y$ that is horizontal with respect
to~$\gamma$.
That is, $p_*\tilde{\xi}=\xi$ and $\gamma(\tilde{\xi})=0$.
For any $\theta\in\mathfrak{u}(1)$, note that $p_*\partial_\theta=0$, which implies that
$p_*[\tilde{\xi},\partial_\theta]=[p_*\tilde{\xi},p_*\partial_\theta]=0$ and
$\gamma([\tilde{\xi},\partial_\theta])=-(p^*\varpi)(\tilde{\xi},\partial_\theta)=0$. Therefore
$[\tilde{\xi},\partial_\theta]=0$ for all~$\theta$.

Associated to~$Y$ is a~complex line bundle~$L$ over~$M$, given by $L=Y\times_{U(1)}\mathbb{C}$.
We write an element of~$L$ as an equivalence class $[y,z]$ with $y\in Y$ and $z\in\mathbb{C}$.
There is a~connection $\nabla$ on~$L$ that is constructed from the connection one-form~$\gamma$ through the following
process.
\begin{itemize}\itemsep=0pt
\item Given any $s\in\Gamma(L)$, def\/ine the map $\tilde{s}:Y\rightarrow\mathbb{C}$ so that $[y,\tilde{s}(y)]=s(p(y))$
for all $y\in Y$.
Then $\tilde{s}$ has the equivariance property
\begin{gather*}
\tilde{s}(y\cdot\lambda)=\lambda^{-1}\tilde{s}(y),\qquad \forall\, y\in Y,\ \lambda\in U(1).
\end{gather*}

\item Conversely, any map $\tilde{s}:Y\rightarrow\mathbb{C}$ with the above equivariance property can be used to
construct a~section~$s$ of~$L$ by setting $s(m)=[y,\tilde{s}(y)]$ for all $m\in M$ and any $y\in Y$ such that $p(y)=m$.

\item Let $\xi\in\mathcal{X}(M)$ be given, and let $\tilde{\xi}$ be its horizontal lift to~$Y$.
If $\tilde{s}:Y\rightarrow\mathbb{C}$ is an equivariant map, then so is $\tilde{\xi}\tilde{s}$.
This follows from the fact that $[\tilde{\xi},\partial_\theta]=0$ for all $\theta\in\mathfrak{u}(1)$.

\item Def\/ine the connection $\nabla$ on~$L$ so that for any $\xi\in\mathcal{X}(M)$ and $s\in\Gamma(L)$, $\nabla_\xi s$
is the section of~$L$ that satisf\/ies
\begin{gather*}
\widetilde{\nabla_\xi s}=\tilde{\xi}\tilde{s}.
\end{gather*}
\end{itemize}

\subsection{The prequantization circle bundle and its inf\/initesimal quantomorphisms}\label{subsec:quanto1}

\begin{Definition}
Let $(M,\omega)$ be a~symplectic manifold.
A~\textit{prequantization circle bundle} for $(M,\omega)$ is a~right principal $U(1)$ bundle
$Y\stackrel{p}{\longrightarrow}M$, together with a~connection one-form~$\gamma$ on~$Y$ satisfying
$d\gamma=\frac{1}{i\hbar}p^*\omega$.
\end{Definition}

\begin{Definition}
Let $(Y_1,\gamma_1)\stackrel{p_1}{\longrightarrow}(M_1,\omega_1)$ and
$(Y_2,\gamma_2)\stackrel{p_2}{\longrightarrow}(M_2,\omega_2)$ be prequantization circle bundles for two symplectic manifolds.
A~dif\/feomorphism $K:Y_1\rightarrow Y_2$ is called a~\textit{quantomorphism} if $K^*\gamma_2=\gamma_1$.
\end{Definition}

Let $K:Y_1\rightarrow Y_2$ be a~quantomorphism.
Notice that for any $\theta\in\mathfrak{u}(1)$, the vector f\/ield $\partial_\theta$ on~$Y_1$ is completely specif\/ied~by
the conditions $\gamma_1(\partial_\theta)=\theta$ and $\partial_\theta\lrcorner d\gamma_1=0$, and the same is true on~$Y_2$.
Since $K^*\gamma_2=\gamma_1$, we see that $K_*\partial_\theta=\partial_\theta$ for all~$\theta$, and so~$K$ is
equivariant with respect to the principal circle actions.

\begin{Definition}
Let $(Y,\gamma)\stackrel{p}{\longrightarrow}(M,\omega)$ be a~prequantization circle bundle.
An \textit{infinitesimal quantomorphism} of $(Y,\gamma)$ is a~vector f\/ield $\zeta\in\mathcal{X}(Y)$ whose f\/low $\phi_t$
on~$Y$ is a~quantomorphism from its domain to its range for each~$t$.
The space of inf\/initesimal quantomorphisms of $(Y,\gamma)$ is denoted by $\mathcal{Q}(Y,\gamma)$.
\end{Definition}

Let $\zeta\in\mathcal{X}(Y)$ have f\/low $\phi_t$.
The connection form~$\gamma$ is preserved by $\phi_t$ if and only if $L_\zeta\gamma=0$.
Therefore the space of inf\/initesimal quantomorphisms of $(Y,\gamma)$ is
\begin{gather*}
\mathcal{Q}(Y,\gamma)=\{\zeta\in\mathcal{X}(Y) \,|\, L_\zeta\gamma=0\}.
\end{gather*}

If $K:Y_1\rightarrow Y_2$ is a~quantomorphism, then it induces a~dif\/feomorphism (in fact, a~symplectomorphism)
$K':M_1\rightarrow M_2$ such that the following diagram commutes:
\begin{gather*}
\xymatrix{Y_1 \ar[r]^K \ar[d]^{p_1} & Y_2 \ar[d]^{p_2} \\
M_1 \ar[r]^{K'} & M_2}
\end{gather*}
This implies that for any $\zeta\in\mathcal{Q}(Y,\gamma)$ with f\/low $\phi_t$, there is a~f\/low $\phi_t'$ on~$M$ that
satisf\/ies $p\circ\phi_t=\phi_t'\circ p$.
If $\zeta'$ is the vector f\/ield on~$M$ with f\/low $\phi_t'$, then $p_*\zeta=\zeta'$.
In other words, elements of $\mathcal{Q}(Y,\gamma)$ descend via $p_*$ to well-def\/ined vector f\/ields on~$M$.

\subsection{The Lie algebra isomorphism}\label{subsec:iso1}

Let $(Y,\gamma)\stackrel{p}{\longrightarrow}(M,\omega)$ be a~prequantization circle bundle.
We will now present an explicit construction of a~Lie algebra isomorphism from $C^\infty(M)$ to $\mathcal{Q}(Y,\gamma)$.
Recall that the vector f\/ield $\partial_{2\pi i}$ on~$Y$ satisf\/ies $\gamma(\partial_{2\pi i})=2\pi i\in\mathfrak{u}(1)$
and~$p_*\partial_{2\pi i}=0$.

\begin{Lemma}
\label{lem:bracketoflift}
For all $f,g\in C^\infty(M)$,
\begin{gather*}
[\tilde{\xi}_f,\tilde{\xi}_g]=\tilde{\xi}_{\{f,g\}}-\frac{1}{2\pi\hbar}p^*\{f,g\}\partial_{2\pi i}.
\end{gather*}
\end{Lemma}

\begin{proof}
It suf\/f\/ices to show that
\begin{gather*}
p_*[\tilde{\xi}_f,\tilde{\xi}_g] =p_*\left(\tilde{\xi}_{\{f,g\}}-\frac{1}{2\pi\hbar}p^*\{f,g\}\partial_{2\pi i}\right)
\end{gather*}
and
\begin{gather*}
\gamma([\tilde{\xi}_f,\tilde{\xi}_g]) =\gamma\left(\tilde{\xi}_{\{f,g\}}-\frac{1}{2\pi\hbar}p^*\{f,g\}\partial_{2\pi i}\right).
\end{gather*}
Using Lemma~\ref{lem:lemma0}, we see that
\begin{gather*}
p_*\left(\tilde{\xi}_{\{f,g\}}-\frac{1}{2\pi\hbar}p^*\{f,g\}\partial_{2\pi i}\right)=\xi_{\{f,g\}}=[\xi_f,\xi_g].
\end{gather*}
Since $p_*\tilde{\xi}_f=\xi_f$ and $p_*\tilde{\xi}_g=\xi_g$, it follows that
$p_*[\tilde{\xi}_f,\tilde{\xi}_g]=[\xi_f,\xi_g]$.
Thus the f\/irst equation is verif\/ied.

Next, note that
\begin{gather*}
\gamma\left(\tilde{\xi}_{\{f,g\}}-\frac{1}{2\pi\hbar}p^*\{f,g\}\partial_{2\pi i}\right)=\frac{1}{i\hbar}p^*\{f,g\},
\end{gather*}
and
\begin{gather*}
\gamma([\tilde{\xi}_f,\tilde{\xi}_g])=-\frac{1}{i\hbar}(p^*\omega)(\tilde{\xi}_f,\tilde{\xi}_g)=\frac{1}{i\hbar}p^*\{f,g\}.
\end{gather*}
Therefore the second equation is also verif\/ied.
\end{proof}

\begin{Lemma}\label{lem:liealghom1}
The map $E:C^\infty(M)\rightarrow\mathcal{X}(Y)$ given~by
\begin{gather*}
E(f)=\tilde{\xi}_f+\frac{1}{2\pi\hbar}p^*f\partial_{2\pi i},\qquad \forall\, f\in C^\infty(M)
\end{gather*}
is a~Lie algebra homomorphism.
\end{Lemma}

\begin{proof}
Let $f,g\in C^\infty(M)$ be arbitrary.
We need to show that
\begin{gather*}
\tilde{\xi}_{\{f,g\}}+\frac{1}{2\pi\hbar}p^*\{f,g\}\partial_{2\pi
i}=\left[\tilde{\xi}_f+\frac{1}{2\pi\hbar}p^*f\partial_{2\pi i},\tilde{\xi}_g+\frac{1}{2\pi\hbar}p^*g\partial_{2\pi
i}\right].
\end{gather*}
Using Lemma~\ref{lem:bracketoflift}, the left-hand side becomes
\begin{gather*}
[\tilde{\xi}_f,\tilde{\xi}_g]+2\frac{1}{2\pi\hbar}p^*\{f,g\}\partial_{2\pi i}.
\end{gather*}
Expanding the right-hand side yields
\begin{gather*}
[\tilde{\xi}_f,\tilde{\xi}_g]+\left[\tilde{\xi}_f,\frac{1}{2\pi\hbar}p^*g\partial_{2\pi
i}\right]+\left[\frac{1}{2\pi\hbar}p^*f\partial_{2\pi
i},\tilde{\xi}_g\right]+\left[\frac{1}{2\pi\hbar}p^*f\partial_{2\pi i},\frac{1}{2\pi\hbar}p^*g\partial_{2\pi i}\right].
\end{gather*}
The fourth term vanishes because $\partial_\theta(p^*f)=\partial_\theta(p^*g)=0$ for any $\theta\in\mathfrak{u}(1)$.
To evaluate the third term, recall that that $[\partial_\theta,\tilde{\xi}]=0$ for any $\theta\in\mathfrak{u}(1)$ and
$\xi\in\mathcal{X}(M)$.
Therefore $[\partial_{2\pi i},\tilde{\xi}_g]=0$, so this term reduces to
\begin{gather*}
-\frac{1}{2\pi\hbar}\big(\tilde{\xi}_g p^*f\big)\partial_{2\pi i}=\frac{1}{2\pi\hbar}p^*\{f,g\}\partial_{2\pi i}.
\end{gather*}
By the same argument, the second term also reduces to
\begin{gather*}
\frac{1}{2\pi\hbar}p^*\{f,g\}\partial_{2\pi i}.
\end{gather*}
Combining these results, we f\/ind that the right-hand side of the desired equation is
\begin{gather*}
\big[\tilde{\xi}_f,\tilde{\xi}_g\big]+2\frac{1}{i\hbar}p^*\{f,g\}\partial_{2\pi i},
\end{gather*}
which equals the left-hand side.
\end{proof}

\begin{Lemma}
\label{lem:target1}
For all $f\in C^\infty(M)$, $E(f)\in\mathcal{Q}(Y,\gamma)$.
\end{Lemma}

\begin{proof}
We need to show that $L_{E(f)}\gamma=0$.
We calculate
\begin{gather*}
L_{E(f)}\gamma=E(f)\lrcorner d\gamma+d(E(f)\lrcorner \gamma)=\frac{1}{i\hbar}p^*(\xi_f\lrcorner
\omega)-\frac{1}{i\hbar}p^*df=0.
\tag*{\qed}
\end{gather*}
\renewcommand{\qed}{}
\end{proof}

So far, we have shown that $E:C^\infty(M)\rightarrow\mathcal{Q}(Y,\gamma)$ is a~Lie algebra homomorphism.
We will now construct a~map $F:\mathcal{Q}(Y,\gamma)\rightarrow C^\infty(M)$, and show that~$E$ and~$F$ are inverses.
This will complete the proof that $C^\infty(M)$ and $\mathcal{Q}(Y,\gamma)$ are isomorphic.

Let $\zeta\in\mathcal{Q}(Y,\gamma)$ be arbitrary.
Then $L_\zeta\gamma=\zeta\lrcorner d\gamma+d(\gamma(\zeta))=0$.
This implies that $\partial_\theta\lrcorner (\zeta\lrcorner d\gamma+d(\gamma(\zeta)))=0$ for any
$\theta\in\mathfrak{u}(1)$.
Since $d\gamma(\zeta,\partial_\theta)=\frac{1}{i\hbar}(p^*\omega)(\zeta,\partial_\theta)=0$, it follows that
$\partial_\theta\lrcorner d(\gamma(\zeta))=L_{\partial_\theta}\gamma(\zeta)=0$.
We can therefore def\/ine the map $F:\mathcal{Q}(Y,\gamma)\rightarrow C^\infty(M)$ so that
\begin{gather*}
-\frac{1}{i\hbar}p^*F(\zeta)=\gamma(\zeta),\qquad \forall\,\zeta\in\mathcal{Q}(Y,\gamma).
\end{gather*}

\begin{Theorem}
\label{thm:liealgiso1}
The map $E:C^\infty(M)\rightarrow\mathcal{Q}(Y,\gamma)$ is a~Lie algebra isomorphism with inverse~$F$.
\end{Theorem}

\begin{proof}
Let $f\in C^{\infty}(M)$ and $\zeta\in\mathcal{Q}(Y,\gamma)$ be arbitrary.
We will show that $F(E(f))=f$ and $E(F(\zeta))=\zeta$.
Using the def\/initions of~$E$ and~$F$, we have
\begin{gather*}
-\frac{1}{i\hbar}p^*F(E(f))=\gamma(E(f))=\gamma\left(\tilde{\xi}_f+\frac{1}{2\pi\hbar}p^*f\partial_{2\pi
i}\right)=-\frac{1}{i\hbar}p^*f.
\end{gather*}
This implies that $F(E(f))=f$.

To show that $E(F(\zeta))=\zeta$, it suf\/f\/ices to show that $\gamma(E(F(\zeta)))=\gamma(\zeta)$ and
$p_*E(F(\zeta))=p_*\zeta$.
By def\/inition,
\begin{gather*}
E(F(\zeta))=\tilde{\xi}_{F(\zeta)}+\frac{1}{2\pi\hbar} p^*F(\zeta)\partial_{2\pi i}=\tilde{\xi}_{F(\zeta)}+\frac{1}{2\pi
i}\gamma(\zeta)\partial_{2\pi i}.
\end{gather*}
It is immediate that $\gamma(E(F(\zeta)))=\gamma(\zeta)$, and that $p_*E(F(\zeta))=\xi_{F(\zeta)}$.
Observe that
\begin{gather*}
\zeta\lrcorner p^*\omega={i\hbar}\zeta\lrcorner d\gamma=-{i\hbar}d(\gamma(\zeta))=p^*(dF(\zeta)),
\end{gather*}
having used $L_\zeta\gamma=0$.
Therefore $(p_*\zeta)\lrcorner \omega=dF(\zeta)$, which implies that $p_*\zeta=\xi_{F(\zeta)}$.
Thus $p_*E(F(\zeta))=p_*\zeta$.
This concludes the proof that $E(F(\zeta))=\zeta$.

Since~$E$ and~$F$ are inverses, and we know from Lemma~\ref{lem:target1} that
$E:C^\infty(M)\rightarrow\mathcal{Q}(Y,\gamma)$ is a~Lie algebra homomorphism, it follows that~$E$ and~$F$ are the
desired Lie algebra isomorphisms.
\end{proof}

The primary goal of Section~\ref{sec:mpc} is to duplicate the above construction for the inf\/initesimal quantomorphisms
of a~metaplectic-c prequantization.
However, before moving on to the metaplectic-c case, we will show how the map~$E$ can be used to represent the elements
of $C^\infty(M)$ as operators on the space of sections of the prequantization line bundle for $(M,\omega)$.
This result will also have an analogue in the metaplectic-c case, which we will discuss in
Section~\ref{subsec:operator2}.

\subsection[An operator representation of $C^\infty(M)$]{An operator representation of $\boldsymbol{C^\infty(M)}$}\label{subsec:operator1}

Let $(L,\nabla)$ be the complex line bundle with connection associated to $(Y,\gamma)$.
One of the goals of the Kostant--Souriau prequantization process is to produce a~representation
$r:C^\infty(M)\rightarrow\End \Gamma(L)$.
To be consistent with quantum mechanics in the case of a~physically realizable system, the map~$r$ is required to
satisfy the following axioms:
\begin{enumerate}[(1)]\itemsep=0pt
\item $r(1)$ is the identity map on $\Gamma(L)$,
\item for all $f,g\in C^\infty(M)$, $[r(f),r(g)]=i\hbar r(\{f,g\})$ (up to sign convention).
\end{enumerate}
These axioms are based on an analysis by Dirac~\cite{dir} of the relationship between classical and quantum mechanical observables.
For more detail in the context of geometric quantization, see, for example, \'Sniatycki~\cite{sn1} or Woodhouse~\cite{w1}.

Recall the association between a~section~$s$ of~$L$ and an equivariant function $\tilde{s}:Y\rightarrow\mathbb{C}$.
We note the following properties.
\begin{itemize}\itemsep=0pt
\item For any $f\in C^\infty(M)$ and $s\in\Gamma(L)$, the equivariant function corresponding to the section~$fs$ is
$\widetilde{fs}=p^*f\tilde{s}$.
\item The vector f\/ield $\partial_{2\pi i}$ has f\/low $R_{\exp(2\pi it)}$.
Thus, for all $y\in Y$,
\begin{gather*}
(\partial_{2\pi i}\tilde{s})(y)=\left.\frac{d{}}{d{t}}\right|_{t=0}\tilde{s}\big(y\cdot e^{2\pi it}\big)=-2\pi i\tilde{s}(y).
\end{gather*}
\end{itemize}
The Kostant--Souriau representation $r:C^\infty(M)\rightarrow\End \Gamma(L)$ is def\/ined by
\begin{gather*}
r(f)s=\left(i\hbar\nabla_{\xi_f}+f\right)s,\qquad \forall\, f\in C^\infty(M),\ s\in\Gamma(L).
\end{gather*}
Using the preceding observations, we see that
\begin{gather*}
\widetilde{r(f)s}=\big(i\hbar\tilde{\xi}_f+p^*f\big)\tilde{s}=\left(i\hbar\tilde{\xi}_f-\frac{1}{2\pi
i}p^*f\partial_{2\pi i}\right)\tilde{s}=i\hbar E(f)\tilde{s}.
\end{gather*}
Since we proved in Lemma~\ref{lem:liealghom1} that $E(\{f,g\})=[E(f),E(g)]$ for all $f,g\in C^\infty(M)$, the following
is immediate.

\begin{Theorem}
The map $r:C^\infty(M)\rightarrow\End \Gamma(L)$ satisfies Dirac axioms $(1)$ and $(2)$.
\end{Theorem}

Thus the same map that provides the isomorphism from $C^\infty(M)$ to $\mathcal{Q}(Y,\gamma)$ also yields the usual
Kostant--Souriau representation of $C^\infty(M)$ as a~space of operators on $\Gamma(L)$.
We will see a~similar result in the case of metaplectic-c prequantization.

\section{Metaplectic-c quantomorphisms}\label{sec:mpc}

Having reviewed the properties of inf\/initesimal quantomorphisms in Kostant--Souriau prequantization, we will now explore
their parallels in metaplectic-c prequantization.
In Sections~\ref{subsec:mpcgp} and~\ref{subsec:mpcquant}, we summarize the prequantization stage of the metaplectic-c
quantization process de\-ve\-loped by Robinson and Rawnsley~\cite{rr1}.
In Section~\ref{subsec:mpcinftl}, we develop our def\/inition for a~metaplectic-c quantomorphism, and use it to def\/ine an
inf\/initesimal metaplectic-c quantomorphism.
The remainder of the paper is dedicated to proving the metaplectic-c analogues of the results presented in
Section~\ref{sec:ks}.

\subsection{The metaplectic-c group}\label{subsec:mpcgp}

Fix a~$2n$-dimensional real vector space~$V$, and equip it with a~symplectic structure~$\Omega$.
Let $\Sp(V)$ be the symplectic group of $(V,\Omega)$; that is, $\Sp(V)$ is the group of linear automorphisms of~$V$ that
preserve~$\Omega$.
The metaplectic group $\Mp(V)$ is the unique connected double cover of $\Sp(V)$.
The metaplectic-c group $\Mp^c(V)$ is def\/ined~by
\begin{gather*}
\Mp^c(V)=\Mp(V)\times_{\mathbb{Z}_2}U(1),
\end{gather*}
where $\mathbb{Z}_2\subset\Mp(V)$ consists of the two preimages of $I\in\Sp(V)$, and $\mathbb{Z}_2\subset U(1)$ is the
usual subgroup $\{1,-1\}$.

Two important group homomorphisms can be def\/ined on $\Mp^c(V)$.
The f\/irst is the projection map $\sigma:\Mp^c(V)\rightarrow\Sp(V)$, which is part of the short exact sequence
\begin{gather*}
1\rightarrow U(1)\rightarrow\Mp^c(V)\stackrel{\sigma}{\longrightarrow}\Sp(V)\rightarrow 1.
\end{gather*}
The second is the determinant map $\eta:\Mp^c(V)\rightarrow U(1)$, which is part of the short exact sequence
\begin{gather*}
1\rightarrow\Mp(V)\rightarrow\Mp^c(V)\stackrel{\eta}{\longrightarrow} U(1)\rightarrow 1.
\end{gather*}
This latter map has the property that if $\lambda\in U(1)\subset\Mp^c(V)$, then $\eta(\lambda)=\lambda^2$.
The Lie algebra $\mathfrak{mp}^c(V)$ is identif\/ied with $\mathfrak{sp}(V)\oplus\mathfrak{u}(1)$ under the map
$\sigma_*\oplus\frac{1}{2}\eta_*$.

\subsection{Metaplectic-c prequantization}\label{subsec:mpcquant}

Let $(M,\omega)$ be a~$2n$-dimensional symplectic manifold.
The symplectic frame bundle for $(M,\omega)$, modeled on $(V,\Omega)$, is denoted $\Sp(M,\omega)$, and is def\/ined
f\/iberwise for each $m\in M$~by
\begin{gather*}
\Sp(M,\omega)_m=\{b:V\rightarrow T_mM \,|\, b \ \text{is an isomorphism and} \ b^*\omega_m=\Omega\}.
\end{gather*}
Then $\Sp(M,\omega)$ is a~right principal $\Sp(V)$ bundle over~$M$, where $g\in\Sp(V)$ acts by precomposition.
Let $\Sp(M,\omega)\stackrel{\rho}{\longrightarrow}M$ be the bundle projection map.

\begin{Definition}
A~\textit{metaplectic-c structure} for $(M,\omega)$ is a~pair $(P,\Sigma)$, where~$P$ is a~right principal $\Mp^c(V)$
bundle $P\stackrel{\Pi}{\longrightarrow}M$, and~$\Sigma$ is a~map $P\stackrel{\Sigma}{\longrightarrow}\Sp(M,\omega)$
that satisf\/ies
\begin{gather*}
\Sigma(q\cdot a)=\Sigma(q)\cdot\sigma(a),\qquad \forall\, q\in P,\ a\in\Mp^c(V),
\end{gather*}
and $\Pi=\rho\circ\Sigma$.
\end{Definition}

We will need the following def\/initions and observations concerning Lie algebras and vector f\/ields.
\begin{itemize}\itemsep=0pt
\item Given $\kappa\in\mathfrak{sp}(V)$, let $\partial_\kappa$ be the vector f\/ield on $\Sp(M,\omega)$ whose f\/low is
$R_{\exp(t\kappa)}$.

\item Given $\alpha\in\mathfrak{mp}^c(V)$, let $\hat{\partial}_\alpha$ be the vector f\/ield on~$P$ whose f\/low is
$R_{\exp(t\alpha)}$.
Under the identif\/ication $\mathfrak{mp}^c(V)=\mathfrak{sp}(V)\oplus\mathfrak{u}(1)$, we can write
$\alpha=\kappa\oplus\tau$ for some $\kappa\in\mathfrak{sp}(V)$ and $\tau\in\mathfrak{u}(1)$.
Naturality of the exponential map and equivariance of~$\Sigma$ with respect to~$\sigma$ ensure that
$\Sigma_*\hat{\partial}_\alpha=\partial_\kappa$.
\end{itemize}

\begin{Definition}
A~\textit{metaplectic-c prequantization} of $(M,\omega)$ is a~triple $(P,\Sigma,\gamma)$, where $(P,\Sigma)$ is
a~metaplectic-c structure for $(M,\omega)$ and~$\gamma$ is a~$\mathfrak{u}(1)$-valued one-form on~$P$ such that:
\begin{enumerate}[(1)]\itemsep=0pt
\item $\gamma$ is invariant under the principal $\Mp^c(V)$ action,
\item $\gamma(\hat{\partial}_\alpha)=\frac{1}{2}\eta_*\alpha$ for all $\alpha\in\mathfrak{mp}^c(V)$,
\item $d\gamma=\frac{1}{i\hbar}\Pi^*\omega$.
\end{enumerate}
\end{Definition}
When~$P$ is viewed as a~bundle over $\Sp(M,\omega)$ with projection map~$\Sigma$, it becomes a~principal circle bundle
with connection one-form~$\gamma$.
The circle that acts on the f\/ibers of~$P$ is the center $U(1)\subset\Mp^c(V)$.

The space of inf\/initesimal quantomorphisms of $(P,\Sigma,\gamma)$ consists of those vector f\/ields on~$P$ whose f\/lows
preserve all of the structures on $(P,\Sigma,\gamma)$.
Note that one of these structures is the map $P\stackrel{\Sigma}{\longrightarrow}\Sp(M,\omega)$, which does not have
a~direct analogue in the Kostant--Souriau case.
We will show how to incorporate this additional piece of information in the next section.

\subsection{Inf\/initesimal metaplectic-c quantomorphisms}\label{subsec:mpcinftl}

As in Section~\ref{subsec:quanto1}, we begin by developing the idea of a~quantomorphism between metaplectic-c
prequantizations.
Let
$(P_1,\Sigma_1,\gamma_1)\stackrel{\Sigma_1}{\longrightarrow}\Sp(M_1,\omega_1)\stackrel{\rho_1}{\longrightarrow}(M_1,\omega_1)$
and $(P_2,\Sigma_2,\gamma_2)\stackrel{\Sigma_2}{\longrightarrow}\Sp(M_2,\omega_2)$
$\stackrel{\rho_2}{\longrightarrow}(M_2,\omega_2)$ be metaplectic-c prequantizations for two symplectic manifolds, and
let $\Pi_j=\rho_j\circ\Sigma_j$ for $j=1,2$.
Let $K:P_1\rightarrow P_2$ be a~dif\/feomorphism.
We will determine the conditions that~$K$ must satisfy in order for it to preserve all of the structures of the
metaplectic-c prequantizations.
First, by analogy with the Kostant--Souriau def\/inition, assume that~$K$ satisf\/ies $K^*\gamma_2=\gamma_1$.

Fix $m\in M_1$, and consider the f\/iber $P_{1m}$.
For any $q\in P_{1m}$, notice that
\begin{gather*}
T_qP_{1m}=\{\xi\in T_qP_1 \,|\, \Pi_{1*}\xi=0\}=\ker d\gamma_{1q}.
\end{gather*}
The same property holds for a~f\/iber of $P_2$ over a~point in $M_2$.
By assumption, $K_*$ is an isomorphism from $\ker d\gamma_{1q}$ to $\ker d\gamma_{2K(q)}$ for all $q\in P_1$.
Therefore $\Pi_2$ is constant on $K(P_{1m})$.
Moreover, since~$K$ is a~dif\/feomorphism, $K(P_{1m})$ is in fact a~f\/iber of $P_2$ over $M_2$, and every f\/iber of $P_2$ is
the image of a~f\/iber of $P_1$.
Thus~$K$ induces a~dif\/feomorphism $K'':M_1\rightarrow M_2$ such that the following diagram commutes:
\begin{gather*}
\xymatrix{P_1 \ar[r]^K \ar[d]^{\Pi_1} & P_2 \ar[d]^{\Pi_2}\\
M_1 \ar[r]^{K''} & M_2
}
\end{gather*}

\begin{Lemma}
\label{lem:symplecto}
The map $K'':M_1\rightarrow M_2$ is a~symplectomorphism.
\end{Lemma}

\begin{proof}
It suf\/f\/ices to show that ${K''}^*\omega_2=\omega_1$.
Using the properties of~$K$, $\gamma_1$ and $\gamma_2$, we calculate
\begin{gather*}
\Pi_1^*({K''}^*\omega_2)=(K''\circ\Pi_1)^*\omega_2=(\Pi_2\circ K)^*\omega_2=K^*(i\hbar d\gamma_2)=i\hbar
d\gamma_1=\Pi_1^*\omega_1.
\end{gather*}
Therefore ${K''}^*\omega_2=\omega_1$, as required.
\end{proof}

Recall from the beginning of Section~\ref{subsec:mpcquant} that an element $b\in\Sp(M_1,\omega_1)_m$ is a~map
$b:V\rightarrow T_mM_1$ such that $b^*\omega_{1m}=\Omega$.
Since~$K''$ is a~symplectomorphism, the composition $K''_*\circ b:V\rightarrow T_{K''(m)}M_2$ satisf\/ies $(K''_*\circ
b)^*\omega_{2K''(m)}=\Omega$, which implies that $K''_*\circ b\in\Sp(M_2,\omega_2)_{K''(m)}$.
Let $\widetilde{K''}:\Sp(M_1,\omega_1)\rightarrow\Sp(M_2,\omega_2)$ be the lift of~$K''$ given~by
\begin{gather*}
\widetilde{K''}(b)=K''_*\circ b,\qquad \forall\, b\in\Sp(M_1,\omega_1).
\end{gather*}
Then $\widetilde{K''}$ is a~dif\/feomorphism, and it is equivariant with respect to the principal $\Sp(V)$ actions.

Thus, if we assume that $K^*\gamma_2=\gamma_1$, we obtain the dif\/feomorphisms $K'':M_1\rightarrow M_2$ and
$\widetilde{K''}:\Sp(M_1,\omega_1)\rightarrow\Sp(M_2,\omega_2)$, where both~$K$ and $\widetilde{K''}$ are lifts of~$K''$.
However,~$K$ is not necessarily a~lift of $\widetilde{K''}$.
Indeed, there might not be any map $K':\Sp(M_1,\omega_1)\rightarrow\Sp(M_2,\omega_2)$ of which~$K$ is a~lift.
A~map~$K$ for which there is no corresponding~$K'$ is constructed in Appendix~\ref{sec:counter},
Example~\ref{ex:nofiber}.
In Section~\ref{subsec:quanto1}, we showed that a~dif\/feomorphism of prequantization circle bundles that preserves the
connection forms must be equivariant with respect to the principal circle actions.
By contrast, Example~\ref{ex:nofiber} demonstrates that it is possible for~$K$ to preserve the prequantization one-forms
without being equivariant with respect to the principal $\Mp^c(V)$ actions.

Suppose we make the additional assumption that $K(q\cdot a)=K(q)\cdot a$ for all $q\in P_1$ and $a\in\Mp^c(V)$.
Then~$K$ induces a~dif\/feomorphism $K':\Sp(M_1,\omega_1)\rightarrow\Sp(M_2,\omega_2)$ that satisf\/ies
$K'\circ\Sigma_1=\Sigma_2\circ K$.
Combining this with the map $K'':M_1\rightarrow M_2$ yields the following commutative diagram:
\begin{gather*}
\xymatrix{P_1 \ar[r]^K \ar[d]^{\Sigma_1} & P_2 \ar[d]^{\Sigma_2}\\
\Sp(M_1,\omega_1) \ar[r]^{K'} \ar[d]^{\rho_1} & \Sp(M_2,\omega_2) \ar[d]^{\rho_2}\\
M_1 \ar[r]^{K''} & M_2
}
\end{gather*}
We now have two maps,~$K'$ and $\widetilde{K''}$, which are dif\/feomorphisms from $\Sp(M_1,\omega_1)$ to
$\Sp(M_2,\omega_2)$.
By construction, $\rho_2\circ K'=\rho_2\circ\widetilde{K''}$, and both~$K'$ and $\widetilde{K''}$ are equivariant with
respect to the principal $\Sp(V)$ actions.
However, it is still possible for~$K'$ and $\widetilde{K''}$ to be dif\/ferent.
A~map~$K$ for which $K'\neq\widetilde{K''}$ is given in Example~\ref{ex:notilde}.

As will be shown in Section~\ref{subsec:iso2}, this potential discrepancy between~$K'$ and $\widetilde{K''}$ must be
prevented in order to construct the desired isomorphism between $C^\infty(M)$ and the inf\/initesimal quantomorphisms.
We therefore propose the following def\/inition.

\begin{Definition}
The dif\/feomorphism $K:P_1\rightarrow P_2$ is a~\textit{metaplectic-c quantomorphism} if
\begin{enumerate}[(1)]\itemsep=0pt
\item $K^*\gamma_2=\gamma_1$,
\item the induced dif\/feomorphism $K'':M_1\rightarrow M_2$ satisf\/ies $\widetilde{K''}\circ\Sigma_1=\Sigma_2\circ K$.
\end{enumerate}
\end{Definition}

Let $K:P_1\rightarrow P_2$ be a~metaplectic-c quantomorphism.
Given our concept of a~quantomorphism as a~dif\/feomorphism that preserves all of the structures of a~metaplectic-c
prequantization, we would expect that~$K$ is equivariant with respect to the $\Mp^c(V)$ actions.
Let $\alpha\in\mathfrak{mp}^c(V)$ be arbitrary, and write $\alpha=\kappa\oplus\tau$ under the identif\/ication of
$\mathfrak{mp}^c(V)$ with $\mathfrak{sp}(V)\oplus\mathfrak{u}(1)$.
The vector f\/ield $\hat{\partial}_\alpha$ on $P_1$ is completely specif\/ied by the conditions
$\gamma_1(\hat{\partial}_\alpha)=\tau$ and $\Sigma_{1*}\hat{\partial}_\alpha=\partial_\kappa$, and the same is true on~$P_2$.
Notice that
\begin{gather*}
\gamma_2(K_*\hat{\partial}_\alpha)=\gamma_1(\hat{\partial}_\alpha)=\tau,
\end{gather*}
and
\begin{gather*}
\Sigma_{2*}K_*\hat{\partial}_\alpha=\widetilde{K''}_*\Sigma_{1*}\hat{\partial}_\alpha=\widetilde{K''}_*\partial_\kappa=\partial_\kappa,
\end{gather*}
where the f\/inal equality follows from the fact that $\widetilde{K''}$ is equivariant with respect to $\Sp(V)$.
Thus $K_*\hat{\partial}_\alpha=\hat{\partial}_\alpha$ for all $\alpha\in\mathfrak{mp}^c(V)$, which implies that~$K$ is
equivariant with respect to $\Mp^c(V)$, as desired.

Now consider a~single metaplectic-c prequantized space
$(P,\Sigma,\gamma)\stackrel{\Sigma}{\longrightarrow}\Sp(M,\omega)\stackrel{\rho}{\longrightarrow}(M,\omega)$ with
$\Pi=\rho\circ\Sigma$.

\begin{Definition}
A~vector f\/ield $\zeta\in\mathcal{X}(P)$ is an \textit{infinitesimal metaplectic-c quantomorphism} if its f\/low $\phi_t$
is a~metaplectic-c quantomorphism from its domain to its range for each~$t$.
\end{Definition}

Let $\zeta\in\mathcal{X}(P)$ have f\/low $\phi_t$.
Property (1) of a~quantomorphism holds for $\phi_t$ if and only if \mbox{$L_\zeta\gamma=0$}.
If we assume that $\phi_t$ satisf\/ies property (1), then we can make the following observations.
\begin{itemize}\itemsep=0pt
\item There is a~f\/low ${\phi_t}''$ on~$M$ such that $\Pi\circ\phi_t={\phi_t}''\circ\Pi$.
The vector f\/ield that it generates on~$M$ is $\Pi_*\zeta$.

\item Lemma~\ref{lem:symplecto} shows that ${\phi_t}''$ is a~family of symplectomorphisms.
Therefore we can lift ${\phi_t}''$ to a~f\/low on $\Sp(M,\omega)$, denoted by $\widetilde{{\phi_t}''}$, where
$\widetilde{{\phi_t}''}(b)=({\phi_t}'')_*\circ b$ for all $b\in\Sp(M,\omega)$.
Let the vector f\/ield on $\Sp(M,\omega)$ that has f\/low $\widetilde{{\phi_t}''}$ be $\widetilde{\Pi_*\zeta}$.

\item Property (2) of a~quantomorphism holds for $\phi_t$ if and only if $\Sigma_*\zeta$ is a~well-def\/ined vector f\/ield
on $\Sp(M,\omega)$ and $\Sigma_*\zeta=\widetilde{\Pi_*\zeta}$.
\end{itemize}
We conclude that the space of inf\/initesimal metaplectic-c quantomorphisms of $(P,\Sigma,\gamma)$ is
\begin{gather*}
\mathcal{Q}(P,\Sigma,\gamma)=\big\{\zeta\in\mathcal{X}(P)\,\big|\, L_\zeta\gamma=0
\
\text{and}
\
\Sigma_*\zeta=\widetilde{\Pi_*\zeta}\big\},
\end{gather*}
where it is understood that the condition $\Sigma_*\zeta=\widetilde{\Pi_*\zeta}$ can only be satisf\/ied if
$\Sigma_*\zeta$ is well def\/ined.
In the next section, we will construct a~Lie algebra isomorphism from $C^\infty(M)$ to $\mathcal{Q}(P,\Sigma,\gamma)$.

\subsection{The Lie algebra isomorphism}\label{subsec:iso2}

We begin with a~procedure, given by Robinson and Rawnsley in~\cite[Section~7]{rr1}, for lifting a~Hamiltonian vector
f\/ield on~$M$ to $\Sp(M,\omega)$ and then to~$P$.
These steps will be used in constructing the isomorphism $E:C^\infty(M)\rightarrow\mathcal{Q}(P,\Sigma,\gamma)$.

Fix $f\in C^\infty(M)$, and let its Hamiltonian vector f\/ield $\xi_f$ have f\/low $\varphi_t$ on~$M$.
We know that~${\varphi_t}_*$ preserves~$\omega$ because $L_{\xi_f}\omega=0$.
Let $\tilde{\varphi}_t$ be the lift of $\varphi_t$ to $\Sp(M,\omega)$ given~by
\begin{gather*}
\tilde{\varphi}_t(b)={\varphi_t}_*\circ b,\qquad \forall\, b\in\Sp(M,\omega),
\end{gather*}
and let the vector f\/ield on $\Sp(M,\omega)$ with f\/low $\tilde{\varphi}_t$ be $\tilde{\xi}_f$.
We have $\rho_*\tilde{\xi}_f=\xi_f$ by construction.
Also, $\tilde{\varphi}_t$ commutes with the right principal $\Sp(V)$ action on $\Sp(M,\omega)$, so
$[\tilde{\xi}_f,\partial_\kappa]=0$ for all $\kappa\in\mathfrak{sp}(V)$.
Now let $\hat{\xi}_f$ be the lift of $\tilde{\xi}_f$ to~$P$ that is horizontal with respect to~$\gamma$.
Then $\Sigma_*\hat{\xi}_f=\tilde{\xi}_f$ and $\gamma(\hat{\xi}_f)=0$.
A~summary of the key properties of $\xi_f$, $\tilde{\xi}_f$ and $\hat{\xi}_f$ is below:
\begin{gather*}
\xymatrix{(P,\Sigma,\gamma) \ar[d]^{\Sigma} & \hat{\xi}_f & \gamma(\hat{\xi}_f)=0,\ \ \Sigma_*\hat{\xi}_f=\tilde{\xi}_f,\ \ \Pi_*\hat{\xi}_f=\xi_f\\
\Sp(M,\omega) \ar[d]^{\rho} & \tilde{\xi}_f &[\tilde{\xi}_f,\partial_\kappa]=0\ \ \forall\,\kappa\in\mathfrak{sp}(V),\ \ \rho_*\tilde{\xi}_f=\xi_f\\
(M,\omega) & \xi_f
}
\end{gather*}

The following is a~consequence of Lemma~\ref{lem:lemma0}.
\begin{Lemma}
For all $f,g\in C^\infty(M)$, $\tilde{\xi}_{\{f,g\}}=[\tilde{\xi}_f,\tilde{\xi}_g]$.
\end{Lemma}

In Section~\ref{subsec:iso1}, we made use of the vector f\/ield $\partial_{2\pi i}$ on~$Y$.
The corresponding object in this context is the vector f\/ield $\hat{\partial}_{2\pi i}$ on~$P$, which satisf\/ies
$\gamma(\hat{\partial}_{2\pi i})=2\pi i$ and $\Sigma_*(\hat{\partial}_{2\pi i})=0$.

\begin{Lemma}
For all $f,g\in C^\infty(M)$,
\begin{gather*}
\big[\hat{\xi}_f,\hat{\xi}_g\big]=\hat{\xi}_{\{f,g\}}-\frac{1}{2\pi\hbar}\Pi^*\{f,g\}\hat{\partial}_{2\pi i}.
\end{gather*}
\end{Lemma}

\begin{proof}
It suf\/f\/ices to show that
\begin{gather*}
\Sigma_*\big[\hat{\xi}_f,\hat{\xi}_g\big] =\Sigma_*\left(\hat{\xi}_{\{f,g\}}-\frac{1}{2\pi\hbar}\Pi^*\{f,g\}\hat{\partial}_{2\pi i}\right)
\end{gather*}
and
\begin{gather*}
\gamma\big(\big[\hat{\xi}_f,\hat{\xi}_g\big]\big) =\gamma\left(\hat{\xi}_{\{f,g\}}-\frac{1}{2\pi\hbar}\Pi^*\{f,g\}\hat{\partial}_{2\pi i}\right).
\end{gather*}
The proof proceeds identically to that of Lemma~\ref{lem:bracketoflift}.
\end{proof}

\begin{Lemma}
\label{lem:liealghom2}
The map $E:C^\infty(M)\rightarrow\mathcal{X}(P)$ given~by
\begin{gather*}
E(f)=\hat{\xi}_f+\frac{1}{2\pi\hbar}\Pi^*f\hat{\partial}_{2\pi i},\qquad \forall\, f\in C^\infty(M)
\end{gather*}
is a~Lie algebra homomorphism.
\end{Lemma}

\begin{proof}
Precisely analogous to Lemma~\ref{lem:liealghom1}.
\end{proof}

\begin{Lemma}
For all $f\in C^\infty(M)$, $E(f)\in\mathcal{Q}(P,\Sigma,\gamma)$.
\end{Lemma}

\begin{proof}
We need to show that $L_{E(f)}\gamma=0$ and $\Sigma_* E(f)=\widetilde{\Pi_*E(f)}$.
The verif\/ication that $L_{E(f)}\gamma=0$ is the same as that in Lemma~\ref{lem:target1}.
Note that $\Sigma_*E(f)=\tilde{\xi}_f$ and $\Pi_*E(f)=\xi_f$, so $\widetilde{\Pi_*E(f)}=\tilde{\xi}_f=\Sigma_* E(f)$.
Thus the necessary conditions are satisf\/ied, and $E(f)\in\mathcal{Q}(P,\Sigma,\gamma)$.
\end{proof}

As before, we will construct an inverse for~$E$, and conclude that~$E$ is a~Lie algebra isomorphism.
Let $\zeta\in\mathcal{Q}(P,\Sigma,\gamma)$ and $\alpha\in\mathfrak{mp}^c(V)$ be arbitrary.
Using an identical argument to the one that precedes Theorem~\ref{thm:liealgiso1}, the fact that
$\hat{\partial}_\alpha\lrcorner (L_\zeta\gamma)=0$ implies that $L_{\hat{\partial}_\alpha}\gamma(\zeta)=0$.
Therefore we can def\/ine $F:\mathcal{Q}(P,\Sigma,\gamma)\rightarrow C^\infty(M)$ so that
\begin{gather*}
-\frac{1}{i\hbar}\Pi^*F(\zeta)=\gamma(\zeta),\qquad \forall\,\zeta\in\mathcal{Q}(P,\Sigma,\gamma).
\end{gather*}

\begin{Theorem}
The map $E:C^\infty(M)\rightarrow\mathcal{Q}(P,\Sigma,\gamma)$ is a~Lie algebra isomorphism with inverse~$F$.
\end{Theorem}

\begin{proof}
Let $f\in C^\infty(M)$ and $\zeta\in\mathcal{Q}(P,\Sigma,\gamma)$ be arbitrary.
From the def\/initions of~$E$ and~$F$, $F(E(f))$ satisf\/ies
\begin{gather*}
-\frac{1}{i\hbar}\Pi^*F(E(f))=\gamma(E(f))=\gamma\left(\hat{\xi}_f+\frac{1}{2\pi\hbar}\Pi^*f\hat{\partial}_{2\pi
i}\right)=-\frac{1}{i\hbar}\Pi^*f.
\end{gather*}
Thus $F(E(f))=f$.

Next, we claim that $\gamma(E(F(\zeta)))=\gamma(\zeta)$ and $\Sigma_* E(F(\zeta))=\Sigma_*\zeta$.
Observe that
\begin{gather*}
E(F(\zeta))=\hat{\xi}_{F(\zeta)}+\frac{1}{2\pi\hbar}\Pi^*F(\zeta)\hat{\partial}_{2\pi
i}=\hat{\xi}_{F(\zeta)}+\frac{1}{2\pi i}\gamma(\zeta)\hat{\partial}_{2\pi i}.
\end{gather*}
It is immediate that $\gamma(E(F(\zeta)))=\gamma(\zeta)$ and $\Sigma_* E(F(\zeta))=\tilde{\xi}_{F(\zeta)}$.
From the def\/inition of $\mathcal{Q}(P,\Sigma,\gamma)$, we know that $\Sigma_*\zeta=\widetilde{\Pi_*\zeta}$.
It remains to show that $\Pi_*\zeta=\xi_{F(\zeta)}$.
We calculate
\begin{gather*}
\zeta\lrcorner \Pi^*\omega=\zeta\lrcorner i\hbar d\gamma=-i\hbar d(\gamma(\zeta))=\Pi^*dF(\zeta).
\end{gather*}
This demonstrates that $(\Pi_*\zeta)\lrcorner \omega=dF(\zeta)$, which implies that $\Pi_*\zeta=\xi_{F(\zeta)}$ as
needed.
Thus we have shown that $E(F(\zeta))=\zeta$, and this completes the proof that~$E$ and~$F$ are inverses.
\end{proof}

If the def\/inition of $\mathcal{Q}(P,\Sigma,\gamma)$ did not include the condition that
$\Sigma_*\zeta=\widetilde{\Pi_*\zeta}$, this proof would fail in the f\/inal step.
We would be able to show that $\Sigma_*E(F(\zeta))=\tilde{\xi}_{F(\zeta)}=\widetilde{\Pi_*\zeta}$, but this vector f\/ield
would not necessarily equal $\Sigma_*\zeta$, and so~$F$ would not be the inverse of~$E$.
This explains why property~(2) of a~metaplectic-c quantomorphism is necessary in order to obtain a~subalgebra of~$\mathcal{X}(P)$ that is isomorphic to $C^\infty(M)$.

\subsection[An operator representation of $C^\infty(M)$]{An operator representation of $\boldsymbol{C^\infty(M)}$}\label{subsec:operator2}

In~\cite{rr1}, Robinson and Rawnsley construct an inf\/inite-dimensional Hilbert space $\mathcal{E}'(V)$ of holomorphic
functions on $V\cong\mathbb{C}^n$, on which the group $\Mp^c(V)$ acts via the metaplectic representation.
They then def\/ine the bundle of symplectic spinors for the prequantized system
$(P,\Sigma,\gamma)\stackrel{\Pi}{\longrightarrow}(M,\omega)$ to be
\begin{gather*}
\mathcal{E}'(P)=P\times_{\Mp^c(V)}\mathcal{E}'(V).
\end{gather*}
We omit the details of the metaplectic representation here; the only fact we need is that the subgroup
$U(1)\subset\Mp^c(V)$ acts on elements of $\mathcal{E}'(V)$ by scalar multiplication.
We write an element of $\mathcal{E}'(P)$ as an equivalence class $[q,\psi]$ with $q\in P$ and $\psi\in\mathcal{E}'(V)$.

Section 7 of \cite{rr1} contains the following construction.
\begin{itemize}\itemsep=0pt
\item Let $s\in\Gamma(\mathcal{E}'(P))$ be given, and def\/ine the map $\tilde{s}:P\rightarrow\mathcal{E}'(V)$ so that
$[q,\tilde{s}(q)]=s(\Pi(q))$ for all $q\in P$.
This map $\tilde{s}$ satisf\/ies the equivariance condition
\begin{gather*}
\tilde{s}(q\cdot a)=a^{-1}\tilde{s}(q),\qquad \forall\, q\in P,\ a\in\Mp^c(V),
\end{gather*}
where the action on the right-hand side is that of the metaplectic representation.

\item Conversely, if $\tilde{s}:P\rightarrow\mathcal{E}'(V)$ is any map with the equivariance property above, it can be
used to def\/ine a~section $s\in\Gamma(\mathcal{E}'(P))$ by setting $s(m)=[q,\tilde{s}(q)]$ for each $m\in M$ and any
$q\in P$ such that $\Pi(q)=m$.

\item Let $f\in C^\infty(M)$ be arbitrary, and recall the lifting $\xi_f\rightarrow\tilde{\xi}_f\rightarrow\hat{\xi}_f$
of $\xi_f$ to~$P$.
A~standard calculation establishes that $[\hat{\xi}_f,\hat{\partial}_\alpha]=0$ for all $\alpha\in\mathfrak{mp}^c(V)$.
Thus, if $\tilde{s}:P\rightarrow\mathcal{E}'(V)$ is an equivariant map, then so is $\hat{\xi}_f\tilde{s}$.

\item Def\/ine the map $D:C^\infty(M)\rightarrow\End \Gamma(\mathcal{E}'(P))$ such that for all $f\in C^\infty(M)$
and $s\in\Gamma(\mathcal{E}'(P))$, $D_fs$ is the section of $\mathcal{E}'(P)$ that satisf\/ies
\begin{gather*}
\widetilde{D_fs}=\hat{\xi}_f\tilde{s}.
\end{gather*}
Further, def\/ine $\delta:C^\infty(M)\rightarrow\End \Gamma(\mathcal{E}'(P))$ by
\begin{gather*}
\delta_fs=D_fs+\frac{1}{i\hbar}fs,\qquad \forall\, f\in C^\infty(M),\ s\in\Gamma(\mathcal{E}'(P)).
\end{gather*}
\end{itemize}
Theorem~7.8 of~\cite{rr1} states that~$\delta$ is a~Lie algebra homomorphism.

We see that the construction of~$D$ precisely parallels the construction of the connection $\nabla$ on the
prequantization line bundle~$L$ associated to a~prequantization circle bundle $(Y,\gamma)$.
As in Section~\ref{subsec:operator1}, we make two observations.
\begin{itemize}\itemsep=0pt
\item For any $s\in\Gamma(\mathcal{E}'(P))$ and $f\in C^\infty(M)$, $\widetilde{fs}=\Pi^*f\tilde{s}$.
\item For any equivariant map $\tilde{s}:P\rightarrow\mathcal{E}'(V)$, $\hat{\partial}_{2\pi i}\tilde{s}=-2\pi
i\tilde{s}$.
\end{itemize}
Therefore
\begin{gather*}
\widetilde{\delta_fs}=\left(\hat{\xi}_f+\frac{1}{2\pi\hbar}\Pi^*f\hat{\partial}_{2\pi i}\right)\tilde{s}=E(f)\tilde{s}.
\end{gather*}
The fact that~$\delta$ is a~Lie algebra homomorphism then follows immediately from Lemma~\ref{lem:liealghom2}.
This construction would apply equally well to any associated bundle where the subgroup $U(1)\subset\Mp^c(V)$ acts on the
f\/iber by scalar multiplication.

\appendix

\section{Example of a~metaplectic-c prequantization}\label{sec:counter}

Recall from Section~\ref{subsec:mpcinftl} that a~metaplectic-c quantomorphism~$K$ between two metaplectic-c
prequantizations
$(P_1,\Sigma_1,\gamma_1)\stackrel{\Sigma_1}{\longrightarrow}\Sp(M_1,\omega_1)\stackrel{\rho_1}{\longrightarrow}(M_1,\omega_1)$
and
$(P_2,\Sigma_2,\gamma_2)\stackrel{\Sigma_2}{\longrightarrow}\Sp(M_2,\omega_2)\stackrel{\rho_2}{\longrightarrow}(M_2,\omega_2)$
is a~dif\/feomorphism $K:P_1\rightarrow P_2$ such that
\begin{enumerate}[(1)]\itemsep=0pt
\item $K^*\gamma_2=\gamma_1$,
\item the induced dif\/feomorphism $K'':M_1\rightarrow M_2$ satisf\/ies $\widetilde{K''}\circ\Sigma_1=\Sigma_2\circ K$,
\end{enumerate}
where $\widetilde{K''}:\Sp(M_1,\omega_1)\rightarrow\Sp(M_2,\omega_2)$ is the lift of~$K''$ given by
\begin{gather*}
\widetilde{K''}(b)=K''_*\circ b,\qquad \forall\, b\in\Sp(M_1,\omega_1).
\end{gather*}
We claimed that condition (1) is insuf\/f\/icient to guarantee that~$K$ is the lift of some map
$K':\Sp(M_1,\omega_1)\rightarrow\Sp(M_2,\omega_2)$.
In particular, a~dif\/feomorphism~$K$ that only satisf\/ies condition (1) might not be equivariant with respect to the
principal $\Mp^c(V)$ actions.
We further claimed that~$K$ might be equivariant and satisfy condition (1), yet fail to satisfy condition (2).
We will now construct examples to support these claims.

Let $M=\mathbb{R}^2\setminus\{(0,0)\}$ with Cartesian coordinates $(p,q)$ and polar coordinates $(r,\theta)$.
Equip~$M$ with the symplectic form $\omega=dp\wedge dq=rdr\wedge d\theta$, and observe that the one-form
$\beta=\frac{1}{2}r^2d\theta$ satisf\/ies $d\beta=\omega$.
Let $V=\mathbb{R}^2$ with basis $\{\hat{x},\hat{y}\}$ and symplectic form $\Omega=\hat{x}^*\wedge\hat{y}^*$, and
consider the global trivialization of the tangent bundle $TM$ such that for all $m\in M$, $T_mM$ is identif\/ied with~$V$
by mapping $\hat{x}\rightarrow\left.\frac{\partial{}}{\partial{p}}\right|_m$ and
$\hat{y}\rightarrow\left.\frac{\partial{}}{\partial{q}}\right|_m$.
Identify $\Sp(M,\omega)$ with $M\times\Sp(V)$ using this trivialization.

Let $P=M\times\Mp^c(V)$, and def\/ine the map $\Sigma:P\rightarrow\Sp(M,\omega)$ by $\Sigma(m,a)=(m,\sigma(a))$ for all
$m\in M$ and $a\in\Mp^c(V)$.
Let $\vartheta_0$ be the trivial connection on the product bundle $M\times\Mp^c(V)$, and let
$\gamma=\frac{1}{i\hbar}\beta+\frac{1}{2}\eta_*\vartheta_0$.
Then $(P,\Sigma,\gamma)$ is a~metaplectic-c prequantization of $(M,\omega)$.
In both of the examples below, we will give a~dif\/feomorphism $K:P\rightarrow P$.

To facilitate the construction in Example~\ref{ex:nofiber}, we introduce a~more explicit representation for elements of
the metaplectic-c group.
By def\/inition, $\Mp^c(V)=\Mp(V)\times_{\mathbb{Z}_2}U(1)$.
The restriction of the projection map $\Mp^c(V)\stackrel{\sigma}{\longrightarrow}\Sp(V)$ to $\Mp(V)$ yields the double
covering $\Mp(V)\stackrel{\sigma}{\longrightarrow}\Sp(V)$.
Write an element of $\Mp^c(V)$ as an equivalence class $[h,e^{2\pi it}]$ with $h\in\Mp(V)$ and $t\in\mathbb{R}$.
In terms of this parametrization, the projection map is given~by
\begin{gather*}
\sigma\big[h,e^{2\pi it}\big]=\sigma(h),
\end{gather*}
and the determinant map $\Mp^c(V)\stackrel{\eta}{\longrightarrow} U(1)$ is given~by
\begin{gather*}
\eta\big[h,e^{2\pi it}\big]=\big(e^{2\pi it}\big)^2.
\end{gather*}

\begin{Example}
\label{ex:nofiber}
We will def\/ine a~dif\/feomorphism $K:P\rightarrow P$ that preserves~$\gamma$, but that does not descend through~$\Sigma$
to a~well-def\/ined map on $\Sp(M,\omega)$.

Let $\mu:\mathbb{R}\rightarrow\Mp(V)$ be any smooth nonconstant path such that $\mu(t+1)=\mu(t)$ for all
$t\in\mathbb{R}$.
Note that the composition $\sigma\circ\mu:\mathbb{R}\rightarrow\Sp(V)$ is also nonconstant.
Now def\/ine $F:\Mp(V)\times U(1)\rightarrow\Mp(V)\times U(1)$~by
\begin{gather*}
F\big(h,e^{2\pi it}\big)=\big(h\mu(2t),e^{2\pi it}\big),\qquad \forall\, h\in\Mp(V),\ t\in\mathbb{R}.
\end{gather*}
This map is a~dif\/feomorphism of $\Mp(V)\times U(1)$, and it descends to a~dif\/feomorphism of $\Mp^c(V)$, which we also
denote~$F$.
For any $[h,e^{2\pi it}]\in\Mp^c(V)$, observe that
\begin{gather*}
\eta\big(F\big[h,e^{2\pi it}\big]\big)=\eta\big[h\mu(2t),e^{2\pi it}\big]=\big(e^{2\pi it}\big)^2=\eta\big[h,e^{2\pi it}\big].
\end{gather*}
This implies that for any $\alpha\in\mathfrak{mp}^c(V)$, $\frac{1}{2}\eta_*F_*\alpha=\frac{1}{2}\eta_*\alpha$.

Def\/ine the dif\/feomorphism $K:P\rightarrow P$ by $K(m,a)=(m,F(a))$ for all $m\in M$ and $a\in\Mp^c(V)$.
Since~$K$ is the identity on~$M$, it preserves~$\beta$.
From the property of~$F$ shown above,~$K$ also preserves $\frac{1}{2}\eta_*\vartheta_0$, and thus it preserves~$\gamma$.
Fix $(m,g)\in\Sp(M,\omega)$, and let $h\in\Mp(V)$ be such that $\sigma(h)=g$.
Then the f\/iber of~$P$ over $(m,g)$ is $P_{(m,g)}=\{(m,[h,e^{2\pi it}]) \,|\, t\in\mathbb{R}\}$.
Notice that
\begin{gather*}
\Sigma\circ K\big(m,\big[h,e^{2\pi it}\big]\big)=\Sigma\big(m,\big[h\mu(2t),e^{2\pi it}\big]\big)=(m,g\sigma(\mu(2t))),
\end{gather*}
which is not constant with respect to~$t$.
Thus $K(P_{(m,g)})$ is not contained within a~single f\/iber of~$P$ over $\Sp(M,\omega)$, which shows that there is no map
$K':\Sp(M,\omega)\rightarrow\Sp(M,\omega)$ such that $K'\circ\Sigma=\Sigma\circ K$.

\end{Example}

If $K:P\rightarrow P$ is equivariant with respect to $\Mp^c(V)$, then it induces a~dif\/feomorphism
$K':\Sp(M,\omega)\rightarrow\Sp(M,\omega)$ that satisf\/ies $K'\circ\Sigma=\Sigma\circ K$.
This map and $\widetilde{K''}$ are both lifts of $K'':M\rightarrow M$, but they might not be the same map.

\begin{Example}
\label{ex:notilde}
We will def\/ine a~dif\/feomorphism $K:P\rightarrow P$ that preserves~$\gamma$ and is equivariant with respect to
$\Mp^c(V)$, but where $K'\neq\widetilde{K''}$.

Let $T_{\lambda}:M\rightarrow M$ be the map that rotates~$M$ about the origin by the angle~$\lambda$, where~$\lambda$ is
not an integer multiple of $2\pi$.
Def\/ine $K:P\rightarrow P$~by
\begin{gather*}
K(m,a)=(T_\lambda(m),a),\qquad \forall\, m\in M,\ a\in\Mp^c(V).
\end{gather*}
Then $K^*\gamma=\gamma$, and $K(q\cdot a)=K(q)\cdot a$ for all $q\in P$ and $a\in\Mp^c(V)$.
The map $K':\Sp(M,\omega)\rightarrow\Sp(M,\omega)$ is given~by
\begin{gather*}
K'(m,g)=(T_\lambda(m),g),\qquad \forall\, m\in M,\ g\in\Sp(V),
\end{gather*}
and the map $K'':M\rightarrow M$ is simply $T_\lambda$.
If we let $T_\lambda$ also denote the automorphism of~$V$ given by rotation about the origin by~$\lambda$, then under
our chosen identif\/ication of $TM$ with $M\times V$, we have
\begin{gather*}
K''_*(m,v)=(T_\lambda(m),T_\lambda(v)),\qquad \forall\, m\in M,\ v\in V.
\end{gather*}
Therefore $\widetilde{K''}:\Sp(M,\omega)\rightarrow\Sp(M,\omega)$ is given~by
\begin{gather*}
\widetilde{K''}(m,g)=(T_\lambda(m),T_\lambda\circ g),\qquad \forall\, m\in M,\ g\in\Sp(V).
\end{gather*}
Hence $K'\neq\widetilde{K''}$.
\end{Example}

\subsection*{Acknowledgements}

The author thanks Alejandro Uribe and Yael Karson for enlightening discussions.
This work was funded in part by an NSERC scholarship.

\pdfbookmark[1]{References}{ref}
\LastPageEnding

\end{document}